\documentclass[aop,seceqn,MSNbibl,citesort,dvips]{arximspdf}
\usepackage{rotating}

\doi{10.1214/10-AOP570}
\volume{39}
\issue{3}
\pubyear{2011}
\firstpage{857}
\lastpage{880}

\makeatletter

\newproclaim{defi}{Definition}[section]
\newtheorem{theorem}[defi]{Theorem}
\newtheorem{cor}[defi]{Corollary}
\newproclaim{rem}[defi]{Remark}
\newtheorem{lem}[defi]{Lemma}
\newtheorem{prop}[defi]{Proposition}

\newcommand{\iint}{\int\!\!\int}
\newcommand{\rmc}{c} 

\newcommand{\R}{\mathbb{R}}
\newcommand{\Rk}{\mathbb{R}^k}
\newcommand{\ent}{\operatorname{Ent}}
\newcommand{\Var}{\operatorname{Var}}
\newcommand{\ep}{\varepsilon}
\newcommand{\T}{\mathbf{T}}

\makeatother

\begin{document}
\begin{frontmatter}

\title{A new characterization of Talagrand's transport-entropy
inequalities and applications}
\runtitle{A new characterization of Talagrand's inequalities}

\begin{aug}
\author[A]{\fnms{Nathael} \snm{Gozlan}\ead
[label=e1]{nathael.gozlan@univ-mlv.fr}},
\author[A]{\fnms{Cyril} \snm{Roberto}\corref{}\thanksref{t1}\ead
[label=e2]{cyril.roberto@univ-mlv.fr}} and
\author[A]{\fnms{Paul-Marie} \snm{Samson}\ead
[label=e3]{paul-marie.samson@univ-mlv.fr}}
\runauthor{N. Gozlan, C. Roberto and P.-M. Samson}
\affiliation{Universit\'e Paris-Est Marne-la-Vall\'ee}
\address[A]{Laboratoire d'Analyse\\
et de Math\'ematiques Appliqu\'ees\\
Universit\'e Paris Est Marne la Vall\'ee\\
(UMR CNRS 8050) 5 bd Descartes\\
77454 Marne la Vall\'ee Cedex 2\\
France\\
\printead{e1}\\
\phantom{E-mail: }\printead*{e2}\\
\phantom{E-mail: }\printead*{e3}} 
\end{aug}

\thankstext{t1}{Supported by the European Research Council
through the ``Advanced Grant'' PTRELSS 228032.}

\received{\smonth{12} \syear{2009}}
\revised{\smonth{6} \syear{2010}}

%
\begin{abstract}
We show that Talagrand's transport inequality is equivalent to a
restricted logarithmic Sobolev inequality. This result clarifies the
links between these two important functional inequalities. As an
application, we give the first proof of the fact that Talagrand's
inequality is stable under bounded perturbations.
\end{abstract}

%

\end{frontmatter}

\section{Introduction}\label{intro}

Talagrand's transport inequality and the logarithmic Sobolev inequality
are known to share important features: they both hold for the Gaussian
measure in any dimension, they enjoy the tensorization property and
they imply Gaussian concentration results. We refer to
\cite{villani,ledoux,ane,gozlan-leonard} for surveys about these notions.
Otto and Villani \cite{otto-villani} proved that the logarithmic
Sobolev inequality implies, in full generality, Talagrand's transport
inequality (see also \cite{bgl}) and under a curvature condition, that
the converse also holds (see also \cite{gozlan}). However, since the
work by Cattiaux and Guillin \cite{cattiaux-guillin}, it is known that the
two inequalities are not equivalent, in general.

In this paper, we prove that Talagrand's transport inequality is
actually equivalent to some restricted form of the logarithmic Sobolev
inequality. Our strategy easily generalizes to other transport
inequalities. As a byproduct, we obtain an elementary and direct proof
of the fact that transport inequalities can be perturbed by bounded functions.

In order to present our main results, we need some definitions and notation.

\subsection{Definitions and notation}\label{sec11}

In all what follows, $c\dvtx\R^k\to\R^+$ is a differentiable function
such that $c(0)=\nabla c(0)=0$. Let $\mu$ and $\nu$ be two
probability measures on
$\R^k$; the \textit{optimal transport cost} between $\nu$ and $\mu$
(with respect to the cost function $c$) is defined by
\[
\mathcal{T}_c(\nu,\mu):=\inf_\pi\biggl\{ \iint c(x-y) \,d\pi(x,y)
\biggr\},
\]
where the infimum is taken over all the probability measures $\pi$ on
$\R^k \times\R^k$ with marginals $\nu$ and $\mu$.
Optimal transport costs are used in a wide class of problems,
in statistics, probability and PDE theory, see \cite{villani}. Here,
we shall focus on the following transport inequality.
\begin{defi}[{[Transportation-cost inequality (\ref{eqiTcC})]}] \label{def:tci}
A probability measure $\mu$ on $\R^k$ satisfies (\ref{eqiTcC}), with
$C>0$, if
{\renewcommand{\theequation}{$\T_c(C)$}
\begin{equation}\label{eqiTcC}\hypertarget{eqiTcClink}
\mathcal{T}_c(\nu,\mu)\leq C H(\nu| \mu)\qquad \forall\nu\in
\mathcal{P}(\R^k),
\end{equation}}

\noindent where
\[
H(\nu|\mu)= \cases{
\displaystyle \int\log\frac{d\nu}{d\mu} \,d\nu, &\quad if $\nu\ll\mu$,\cr
+\infty, &\quad otherwise,}
\]
is the relative entropy of $\nu$ with respect to $\mu$ and $ \mathcal
{P}(\R^k)$ is the set of all probability measures on $\Rk$.
\end{defi}

The inequality (\ref{eqiTcC}) implies concentration results as shown by
Marton \cite{marton}, see also \cite{bobkov-gotze,ledoux}
and \cite{gozlan-leonard} for a full introduction to
this notion.

The quadratic cost $c(x)=|x|^2/2$ (where \mbox{$|\cdot|$} stands for the
Euclidean norm) plays a special role.
In this case, we write \hyperlink{eqiTcClink}{($\T_2(C)$)} and say that Talagrand's transport,
or the quadratic transport, inequality is satisfied. Talagrand proved
in \cite{talagrand}, among other results, that the standard Gaussian
measure satisfies \hyperlink{eqiTcClink}{($\T_2(1)$)} in all dimensions. In turn, inequality
\hyperlink{eqiTcClink}{($\T_2(C)$)} implies dimension free Gaussian concentration results. Recently,
the first author showed that the converse is also true, namely that a
dimension free Gaussian concentration result implies \hyperlink{eqiTcClink}{($\T_2(C)$)} \cite{gozlan}.

Now, we introduce the notion of restricted logarithmic Sobolev
inequalities. To that purpose, we need first to
define $K$-semi-convex functions.
\begin{defi}[($K$-semi-convex function)]
A function $f\dvtx\R^k\to\R$ is $K$-semi-convex ($K \in\R$) for the
cost function $c$ if for all
$\lambda\in[0,1]$, and all $x,y\in\R^k$
%
%
\setcounter{equation}{2}
\begin{eqnarray}\label{K semi-convex cost c 1}
f\bigl(\lambda x+(1-\lambda)y\bigr)&\leq&\lambda f(x)+(1-\lambda)f(y)+\lambda
Kc\bigl((1-\lambda)(y-x)\bigr)\nonumber\\[-8pt]\\[-8pt]
&&{}+(1-\lambda) Kc\bigl(\lambda(y-x)\bigr).\nonumber
\end{eqnarray}
\end{defi}

As shown in Proposition \ref{prop:sc} below, for differentiable
functions, (\ref{K semi-convex cost c 1})
is equivalent to the condition
\[
f(y)\geq f(x)+\nabla f(x)\cdot(y-x)-Kc(y-x)\qquad \forall x,y\in\R^k.
\]
The reader\vspace*{1pt} might see the semi-convexity as an answer to the question:
how far is the function $f$ from being convex?
The quadratic case $c(x)=\frac{1}{2}|x|^2$ is particularly
enlightening since a function $f$ is $K$-semi-convex if
and only if $x \mapsto f(x)+\frac{K}{2}|x|^2$ is convex. Note that the
semi-convexity can be related to the notion of convexity-defect, see,
for example, \cite{barthe-kolesnikov} and references therein where it
is largely discussed and used.
Note also that our definition differs from others, such as \cite
{villani}, Definition~10.10, or
\cite{evans}, Lemma 3 in Chapter 3, page 130.

Dealing only with semi-convex functions leads to the following definition.
\begin{defi}[{[Restricted (modified) logarithmic Sobolev inequality]}]
A probability measure $\mu$ on $\R^k$ verifies \textit{the restricted
logarithmic Sobolev inequality} with constant $C>0$, in short
(\ref{eqrLSIC}), if for all $0\leq K<\frac{1}{C}$ and all
$K$-semi-convex $f\dvtx\R^k\to\R$,
{\renewcommand{\theequation}{$\mathbf{rLSI}(C)$}
\begin{equation}\label{eqrLSIC}
\ent_{\mu}(e^f)\leq\frac{2C}{(1-KC)^2}
\int|\nabla f|^2 e^f \,d\mu,
\end{equation}}

\noindent
where $\ent_{\mu}(g):=\int g\log g \,d\mu-\int g \,d\mu\log\int g
\,d\mu$.
More generally, a probability measure $\mu$ on $\R^k$ verifies the
\textit{restricted modified logarithmic Sobolev inequality}
with constant $C>0$ for the cost $c$, in short (\ref{eqrMLSIcC}),
if for all $K\geq0$, $\eta>0$ with $ \eta+K<1/C$ and all $K$-semi-convex
$f\dvtx\R^k\to\R$ for the cost $c$,
{\renewcommand{\theequation}{$\mathbf{rMLSI}(c,C)$}
\begin{equation}\label{eqrMLSIcC}\hypertarget{eqrMLSIcClink}\quad
\ent_{\mu}(e^f)\leq\frac{\eta}{1-C(\eta+ K)} \int
c^*\biggl(\frac{\nabla f}{\eta}\biggr) e^f \,d\mu,
\end{equation}}

\noindent
where $c^*(u):=\sup_{h\in\R^k} \{ u\cdot h-c(h) \}$ and
$u\cdot h$ is the usual scalar product in $\Rk$.
\end{defi}

Note that (\ref{eqrMLSIcC}) reduces to (\ref{eqrLSIC}) for
$c(x)=c^*(x)=\frac{1}{2} |x|^2$, optimizing over $\eta$.

Without the restriction on the set of $K$-semi-convex functions, the
first inequality
corresponds to the usual logarithmic Sobolev inequality introduced by
Gross \cite{gross} (see also \cite{stam}).
For the second one (without the restriction), we recognize the modified
logarithmic Sobolev inequalities
introduced first by Bobkov and Ledoux~\cite{bobkov-ledoux}, with
$c^*(t)=2|t|^2/(1-\gamma)$ for $|t| \leq\gamma$ and $c^*(t)=+\infty
$ otherwise, $t\in\R$, in order to recover the celebrated result by
Talagrand \cite{talagrand91} on the concentration phenomenon for
products of exponential measures.
Gentil, Guillin and Miclo \cite{ggm} established modified logarithmic
Sobolev inequalities for products of the probability measures
$d\nu_p(t)=e^{-|t|^p}/Z_p$, $t\in\R$ and $p \in(1,2)$, with
$c^*(t)$ that compares to $\max(t^2,|t|^q)$ where $q=p/(p-1) \in
(2,\infty)$
is the dual exponent of $p$. In a subsequent paper \cite{ggm2}, they
generalized their results to a large class of measures with tails
between exponential and Gaussian (see also \cite{barthe-roberto} and
\cite{gozlan07}). In \cite{ggm}, the authors also prove that the
modified logarithmic Sobolev inequality [without the restriction, and
with $c^*(t)$ that compares to $\max(t^2,|t|^q)$] implies the
corresponding transport inequality (\ref{eqiTcC}).

Our results below show that the functional inequalities \hyperlink{eqrMLSIcClink}{($\mathbf
{rMLSI}(c, \cdot)$)} and \hyperlink{eqiTcClink}{($\T_c( \cdot)$)} are equivalent (up to
universal factors in the constants). To give a more complete
description of this equivalence, let us consider yet another type of
logarithmic Sobolev inequalities that we call inf-convolution
logarithmic Sobolev inequality.
\begin{defi}[(Inf-convolution logarithmic Sobolev inequality)]
A probability measure $\mu$ on $\R^k$ verifies \textit{the
inf-convolution logarithmic Sobolev inequality} with constant $C>0$, in
short (\ref{eqICLSIcC}), if for all $\lambda\in(0,1/C)$ and all
$f\dvtx\R^k\to\R$,
{\renewcommand{\theequation}{$\mathbf{ICLSI}(c,C)$}
\begin{equation}\label{eqICLSIcC}\hypertarget{eqICLSIcClink}
\ent_{\mu}(e^f)\leq\frac{1}{1-\lambda C} \int(f-Q^\lambda
f) e^f \,d\mu,
\end{equation}}

\noindent where $Q^\lambda f\dvtx\R^k \rightarrow\R$ denotes the
infimum-convolution of $f$:
\[
Q^\lambda f(x) =\inf_{y\in R^k} \{f(y)+\lambda c(x-y)\}.
\]
\end{defi}
%


\subsection{Main results}\label{sec12}

Our first main result is the following.
%
%
\begin{theorem}\label{main-result2}
Let $\alpha\dvtx\R\to\R^+$ be a convex symmetric function of class
$C^1$ such that $\alpha(0)=\alpha'(0)=0$,
$\alpha'$ is concave on $\R^+$. 
Define $c(x)=\sum_{i=1}^k \alpha(x_i)$ and let
$\mu$ be a probability measure on $\R^k$. The following propositions
are equivalent:
\begin{enumerate}[(3)]
\item[(1)] There exists $C_1>0$ such that $\mu$ verifies the inequality
\hyperlink{eqiTcClink}{($\T_c(C_1)$)}.
\item[(2)] There exists $C_2>0$ such that $\mu$ verifies the inequality
\hyperlink{eqICLSIcClink}{($\mathbf{ICLSI}(c,C_2)$)}.
\item[(3)] There exists $C_3>0$ such that $\mu$ verifies the inequality
\hyperlink{eqrMLSIcClink}{($\mathbf{rMLSI}(c,C_3)$)}.
\end{enumerate}
The constants $C_1$, $C_2$ and $C_3$ are related in the following way:
\begin{eqnarray*}
(1)\Rightarrow(2)\Rightarrow(3)  \qquad\mbox{with } C_1&=&C_2=C_3, \\
(3)\Rightarrow(1)  \qquad\mbox{with } C_1&=&8C_3.
\end{eqnarray*}
%
\end{theorem}

The typical example of function $\alpha$ satisfying the setting of
Theorem \ref{main-result2}
is a smooth version of $\alpha(x) = \min(x^2,x^p)$, with $p \in[1,2]$.
%


The first part $(1) \Rightarrow(2)\Rightarrow(3)$ actually holds in a
more general setting (see Theorem \ref{th:main1}), it is proven in
Section \ref{sec:12}. Moreover, the inequality (\ref{eqICLSIcC})
has a meaning even if $\R^k$ is replaced by an abstract metric space
$X$.
The proof of the second part $(3)\Rightarrow(1)$ is given in Section
\ref{sec:212}.
It uses the Hamilton--Jacobi approach of \cite{bgl} based on explicit
computations on the sup-convolution semi-group (Hopf--Lax formula).
An alternative proof of $(3)\Rightarrow(1)$, with a worst constant,
is given in the subsequent Section \ref{sec:alternative} in the
particular case of
the quadratic cost $c(x)= |x|^2/2$. We believe that such an approach
may lead to further developments in the future and so that it is worth
mentioning it.

In order to keep the arguments as clean as possible and to go straight
to the proofs, we decided to collect most of results
on semi-convex functions, and most of the technical lemmas, in an
independent section (Section \ref{sec5}).

Finally, we present some extensions and comments in Section \ref{sec6}.
We first
give an extension of our main Theorem \ref{th:main1} to Riemannian
manifolds verifying a certain curvature condition (see Theorem \ref
{extensionthm}). Then, in Section \ref{autres log-sob}, we show that
other types of\vadjust{\goodbreak} logarithmic Sobolev inequalities can be derived from
transport inequalities (see Theorem \ref{th:main1bis}). The last
Section \ref{poincare} is a discussion on the links between Poincar\'e
inequality and (restricted) modified logarithmic Sobolev inequality.

Let us end this Introduction with an important application of Theorem
\ref{main-result2}.
It is well known that many functional inequalities of Sobolev type are
stable under bounded perturbations. The first perturbation property of
this type was established by Holley and Stroock in \cite{HS87} for the
logarithmic Sobolev inequality.
\begin{theorem}[(Holley--Stroock)]\label{HS}
Let $\mu$ be a probability measure verifying the logarithmic Sobolev
inequality with a constant $C>0$ [$\mathbf{LSI}(C)$ for short]:
\[
\mathrm{Ent}_\mu(f^2)\leq C \int|\nabla f|^2 \,d\mu\qquad\forall f.
\]
Let $\varphi$ be a bounded function; then the probability measure
$d\tilde{\mu}=\frac{1}{Z}e^{\varphi} \,d\mu$ verifies $\mathbf
{LSI}$ with the constant $\tilde{C}=e^{\mathrm{Osc}(\varphi)}C$,
where the oscillation of $\varphi$ is defined by
\[
\mathrm{Osc}(\varphi)=\sup\varphi- \inf\varphi.
\]
\end{theorem}

A longstanding open question was to establish such a property for
transport inequalities. We have even learned from Villani that this
question was one of the initial motivations behind the celebrated work
\cite{otto-villani}. The representation furnished by Theorem \ref
{main-result2} is the key that enables us to give the first bounded
perturbation property for transport inequalities. The following
corollary is our second main result.
\begin{cor}
Let $\alpha$ be a convex symmetric function of class $C^1$ such that
$\alpha(0)=\alpha'(0)=0$,
$\alpha'$ is concave on $\R^+$. Let $c(x)=\sum_{i=1}^k \alpha(x_i)$ and
$\mu$ be a probability measure on $\R^k$.
Assume that $\mu$ verifies (\ref{eqiTcC}). Let $\varphi\dvtx\R^k\to\R$ be
bounded and define $d\tilde{\mu}(x)=\frac{1}{Z}e^{\varphi(x)} \,d\mu
(x)$, where $Z$ is the normalization constant. Then $\tilde{\mu}$
verifies \hyperlink{eqiTcClink}{($\T_c(8 C e^{\mathrm{Osc}(\varphi)})$)} where
$\mathrm{Osc} (\varphi)= \sup\varphi- \inf\varphi$.
\end{cor}
\begin{pf}
The proof below is a straightforward adaptation of the original proof
of Theorem \ref{HS}.
Using the following representation of the entropy
\[
\ent_\mu(g ) = \inf_{t >0} \biggl\{ \int\biggl(g \log
\biggl( \frac{g}{t} \biggr) - g + t \biggr) \,d\mu\biggr\}
\]
with $g=e^f$, we see that [since $g \log( \frac{g}{t} )
- g + t \geq0$]
\[
\ent_{\tilde{\mu}} (g) \leq\frac{e^{\sup\varphi
}}{Z} \ent_{\mu} (g) .
\]
From the first part of Theorem \ref{main-result2}, it follows that for
all $K\geq0$, $\eta>0$, with $\eta+K<1/C$ and all $K$-semi-convex
functions $f$ for the cost $c$,
\begin{eqnarray*}
\ent_{\tilde{\mu}} (e^f)
& \leq &
\frac{e^{\sup\varphi}}{Z} \frac{\eta}{1-C(\eta+ K)} \int c^*
\biggl(\frac{\nabla f}{\eta}\biggr) e^f \,d\mu\\
& \leq &
\frac{\eta e^{{\mathrm{Osc}}(\varphi)}}{1-C(\eta+ K)} \int c^*
\biggl(\frac{\nabla f}{\eta}\biggr) e^f \,d\tilde{\mu} .
\end{eqnarray*}
Let $u=e^{{\mathrm{Osc}}(\varphi)}$ and $c_u(x):=uc(x/u)$, $x\in\Rk
$. Let $f$ be a $K$-semi-convex function for the cost $c_u$. Since
$u\geq1$ the convexity of $\alpha$ yields $c_u(x)\leq c(x)$, $x\in
\Rk$. Hence, $f$ is a $K$-semi-convex function for the cost $c$.
Observing that $c^*_u(x)= u c^*(x), x\in\Rk$, from the above
inequality, it follows that $\tilde{\mu}$ verifies the inequality
\hyperlink{eqrMLSIcClink}{($\mathbf{rMLSI}(c_u,C)$)}. Then, the second part of Theorem \ref
{main-result2} implies that $\tilde{\mu}$ verifies
\hyperlink{eqiTcClink}{($\T_{c_u}(8C)$)}. From point (i) of the technical Lemma \ref
{lem:tec2}, one has $uc(x/u)\geq c(x)/u$ for $u\geq1$, $x\in\Rk$.
This inequality completes the proof.
\end{pf}
\begin{rem}
After the preparation of this work, we have learned from E.~Milman that
he has obtained in \cite{Mil09e} new perturbation results for various
functional inequalities on a Riemannian manifold equipped with a
probability measure $\mu$ absolutely continuous with respect to the
volume element. His results also cover transport inequalities but are
only true under an additional curvature assumption. To be more precise,
suppose that $\mu$ verifies say \hyperlink{eqiTcClink}{($\T_2(C)$)} and consider another
probability measure of the form $d\tilde{\mu}(x)=e^{-V(x)} \,dx$ such
that
\[
\mathrm{Ric} + \operatorname{Hess} V\geq-\kappa,
\]
for some $\kappa\geq0$. Then if $C>\frac{\kappa}{2}$ and if $\mu$
and $\tilde{\mu}$ are close in some sense to each other, then $\tilde
{\mu}$ verifies \hyperlink{eqiTcClink}{($\T_{2}(\tilde{C})$)} for some $\tilde{C}$ depending
only on $C$, $\kappa$ and on the ``distance'' between $\mu$ and
$\tilde{\mu}$.
Actually, the curvature assumption above makes possible to go beyond
the classical Holley--Stroock property and to work with measures $\tilde
{\mu}$ which are more serious perturbations of $\mu$. Proofs of these
results are based on the remarkable equivalence between concentration
and isoperimetric inequalities under curvature bounded from below,
discovered by Milman in \cite{Mil09d}.
\end{rem}

\section{From transport inequalities to restricted modified
logarithmic Sobolev inequalities} \label{sec:12}

In this section, we prove the first part $(1) \Rightarrow
(2)\Rightarrow(3)$ of Theorem~\ref{main-result2}. As mentioned in the
\hyperref[intro]{Introduction},
this implication holds in a more general setting as we explain now.

Let $X$ denote a Polish space equipped with the Borel $\sigma
$-algebra. Then the optimal transport cost between two probability
measures $\mu$ and $\nu$ on $X$, with cost $\rmc \dvtx X \times X
\to
\mathbb{R}^+$ is
\[
\mathcal{T}_{\rmc}(\nu,\mu):=\inf_\pi\iint\rmc(x,y) \,d\pi(x,y),
\]
where the infimum is taken over all probability measures $\pi$ on $X
\times X$ with marginals $\nu$ and $\mu$. Assume $\rmc$ is
symmetric so that $\mathcal{T}_{\rmc}(\nu,\mu)=\mathcal
{T}_{\rmc}(\mu,\nu)$.
The transport inequality \hyperlink{eqiTcClink}{($\T_{\rmc}(C)$)} is defined accordingly as in
Definition \ref{def:tci}.
For $f\dvtx X \to\mathbb{R}$ and $\lambda>0$, the inf-convolution
$Q^\lambda f \dvtx X \to\mathbb{R}$ is given by
\[
Q^\lambda f(x) = \inf_{y\in X} \{ f(y) + \lambda{\rmc} (
x,y) \}.
\]
The first part of Theorem \ref{main-result2} will be a consequence of
the following general result.
\begin{theorem} \label{th:main1}
Let ${\rmc}\dvtx X \times X \rightarrow\R^+$ be a symmetric
continuous function.
Let $\mu$ be a probability measure on $X$ satisfying \hyperlink{eqiTcClink}{($\T_{\rmc}(C)$)}
for some $C>0$.
Then for all functions $f \dvtx X \to\mathbb{R}$ and all $\lambda\in
(0,1/C)$, it holds
\[
\ent_\mu(e^f) \leq\frac{1}{1 - \lambda C} \int
(f - Q^\lambda f) e^f \,d\mu.
\]
Assume moreover that ${\rmc}(x,y)=c(x-y)$, $x,y\in\R^k$, where
$c\dvtx\R
^k\rightarrow\R^+ $ is a differentiable function such that
$c(0)=\nabla c(0)=0$. Then $\mu$ verifies the inequality \hyperlink{eqrMLSIcClink}{($\mathbf
{rMLSI}(c,C)$)}.
\end{theorem}
\begin{pf}
Fix $f\dvtx X \to\mathbb{R}$, $\lambda\in(0,1/C)$,
and define $d\nu_{f}=\frac{e^f}{\int e^f \,d\mu} \,d\mu$. One has
\begin{eqnarray*}
H(\nu_{f}|\mu)
& = &
\int\log\biggl( \frac{e^f}{\int e^f \,d\mu}\biggr)\frac{e^f}{\int
e^f \,d\mu} \,d\mu
=
\int f \,d\nu_{f}-\log\int e^f \,d\mu\\
&
\leq&\int f \,d\nu_{f}-\int f \,d\mu,
\end{eqnarray*}
where the last inequality comes from Jensen inequality.
Consequently, if $\pi$ is a probability measure on $X\times X$ with
marginals $\nu_f$ and $\mu$
\[
H(\nu_{f}| \mu)\leq\iint\bigl(f(x)-f(y)\bigr) \,d\pi(x,y).
\]
It follows from the definition of the inf-convolution function that
$f(x)-f(y) \leq f(x) -Q^\lambda f(x) + \lambda{\rmc}(x,y)$, for all
$x,y \in X$. Hence,
\[
H(\nu_{f}| \mu)
\leq\iint\bigl(f(x) - Q^\lambda f(x) \bigr) \,d\pi(x,y) + \lambda
\iint{\rmc}(x,y) \,d\pi(x,y),
\]
and optimizing over all $\pi$ with marginals $\nu_f$ and $\mu$
\begin{eqnarray*}
H(\nu_{f}| \mu)&=& \int(f - Q^\lambda f) \,d\nu_f +
\lambda\mathcal{T}_c(\nu_f,\mu) \\
&
\leq&\frac{1}{\int e^f \,d\mu} \,d\mu\int(f - Q^\lambda f
) e^f \,d\mu+ \lambda C H(\nu_{f}| \mu) .
\end{eqnarray*}
The first part of Theorem \ref{th:main1} follows by noticing that
$(\int e^f \,d\mu) H(\nu_{f}| \mu) = \ent_\mu
(e^f)$.
Then the proof of Theorem \ref{th:main1} is completed by applying
Lemma \ref{lem:easy} below.
\end{pf}
%
%
%
%
\begin{lem} \label{lem:easy}
Let $c\dvtx\mathbb{R}^k \to\mathbb{R}^+$ be a differentiable function
such that $c(0)=\nabla c(0)=0$ and define $c^*(x)=\sup_y \{ x \cdot y
- c(y) \}\in\R\cup\{+\infty\}$, $x\in\R^k$.
Then, for any $K$-semi-convex differentiable function $f \dvtx\mathbb
{R}^k \to\mathbb{R}$ for the cost $c$, it holds
\[
f(x) - Q^{K+\eta}f(x) \leq\eta c^* \biggl( - \frac{\nabla f(x)
}{\eta} \biggr)\qquad
\forall x \in\mathbb{R}^k, \forall\eta>0 .
\]
\end{lem}
\begin{pf}
Fix a $K$-semi-convex differentiable function $f \dvtx\mathbb{R}^k \to
\R
$. Also fix $x \in\mathbb{R}^k$ and $\eta>0$.
By Proposition \ref{prop:sc} and\vadjust{\goodbreak} the Young inequality $X\cdot Y \leq
\eta c^*(\frac{X}{\eta}) + \eta c(Y)$, we have
\[
f(x) - f(y) - Kc(y-x)
\leq
- \nabla f(x) \cdot(y-x)
\leq
\eta c^* \biggl( - \frac{\nabla f (x) }{\eta} \biggr) + \eta c( y-x ).
\]
Hence, for any $y \in\mathbb{R}^k$,
\[
f(x) - f(y) -(K+\eta)c(y-x) \leq\eta c^* \biggl( - \frac{\nabla f
(x) }{\eta} \biggr) .
\]
This yields the expected result.
\end{pf}

\section{From restricted modified logarithmic Sobolev inequalities to
transport inequalities---I: Hamilton--Jacobi approach} \label{sec:212}

In this section, we prove the second part $(3) \Rightarrow(1)$ of
Theorem \ref{main-result2}.
The proof is based on the approach of Bobkov,
Gentil and Ledoux \cite{bgl}, using the Hamilton--Jacobi equation. We
will use the following notation: given a convex function $\alpha\dvtx\R
\rightarrow\R^+$ with $\alpha(u)\neq0$ for $u\neq0$, we define
%
%
\begin{equation} \label{eq:omega}
\omega_\alpha(x) = \sup_{u >0} \frac{\alpha(ux)}{\alpha(u)}\qquad
\forall x\in\R.
\end{equation}
\begin{pf*}{Proof of $(3) \Rightarrow(1)$ of Theorem
\ref{main-result2}}
Let $f \dvtx\mathbb{R}^k \to\mathbb{R}$ be a bounded continuous function.
For $x \in\mathbb{R}^k$ and $t \in(0,1)$, define
\[
P_tf(x) = \sup_{y \in\mathbb{R}^k} \biggl\{ f(y) - t c \biggl( \frac
{x-y}{t} \biggr) \biggr\} .
\]
It is well known that $u_{t}=P_{t}f$ verifies the following
Hamilton--Jacobi equation (see, e.g., \cite{evans}): for almost every
$x\in\R^k$ and almost every $t \in(0,+\infty)$,
\[
\cases{
\partial_t u_t (x) = c^* ( -\nabla u_t(x) ), \cr
u_0 = f.}
\]
To avoid lengthy technical arguments, we assume in the sequel that
$P_{t}f$ is continuously differentiable in space and time and that the
equation above holds for all $t$ and $x$. We refer to
\cite{LV07}, proof of Theorem 1.8, or \cite{villani}, proof of Theorem
22.17, for a
complete treatment of the problems arising from the nonsmoothness of $P_{t}f$.
Defining $Z(t)=\int e^{\ell(t )P_{1-t} f} \,d\mu$, where $\ell$ is a
smooth nonnegative function on $\R^+$ with $\ell(0)=0$ that will be
chosen later, one gets
\begin{eqnarray*}
Z'(t)
& = &
\int\biggl( \ell'(t)P_{1-t} f + \ell(t) \,\frac{\partial}{\partial
t}P_{1-t} f \biggr) e^{\ell(t ) P_{1-t} f} \,d\mu\\
& = &
\int\ell'(t) P_{1-t} f e^{\ell(t) P_{1-t} f} \,d\mu- \ell(t) \int
c^* ( {\nabla P_{1-t} f} )e^{\ell(t) P_{1-t} f} \,d\mu.
\end{eqnarray*}
On the other hand,
\[
\ent_\mu\bigl( e^{\ell(t) P_{1-t} f} \bigr) = \ell(t) \int
P_{1-t} f e^{\ell(t) P_{1-t} f} \,d\mu- Z(t) \log Z(t) .
\]
Therefore provided $\ell'(t)\neq0$,
%
%
\begin{eqnarray} \label{eq:step1}
\ent_\mu\bigl( e^{\ell(t) P_{1-t} f} \bigr)
& = &
\frac{\ell(t)}{\ell'(t)} Z'(t) -Z(t) \log Z(t) \nonumber\\[-8pt]\\[-8pt]
&&{}
+ \frac{\ell(t)^2}{\ell'(t)} \int c^* ( {\nabla P_{1-t} f}
)e^{\ell(t) P_{1-t} f} \,d\mu.\nonumber
\end{eqnarray}
By Lemma \ref{lem:semiconv} [with $A=\ell(t)(1-t)$ and $B=1-t$], the
function $g=\ell(t) P_{1-t}f$
is $K(t)$ semi-convex for the cost function $c(x)=\sum_{i=1}^k\alpha
(x_i)$, $x\in\R^k$, where $K(t)=4 \ell(t)(1-t) \omega_\alpha
( \frac{1}{2(1-t)} )$.
Hence, we can apply the restricted logarithmic Sobolev inequality to
get that for any $\eta>0$, any $t \in(0,1)$ such that
$K(t) + \eta<1/C_3$,\setcounter{footnote}{1}\footnote{Note that this condition is not empty
since $K(0)=0$.}
\begin{eqnarray*}
\ent_\mu\bigl( e^{\ell(t) P_{1-t} f} \bigr)
&
\leq&\frac{\eta}{1-(K(t) + \eta)C_3} \int c^* \biggl( \frac{\ell
(t) \nabla P_{1-t} f}{\eta} \biggr) e^{\ell(t) P_{1-t} f} \,d\mu\\
&
\leq&
\frac{\eta\omega_{\alpha^*} ( {\ell(t)}/{\eta}
)}{1-(K(t) + \eta)C_3} \int c^* ( {\nabla P_{1-t} f} )
e^{\ell(t) P_{1-t} f} \,d\mu,
\end{eqnarray*}
since $c^*(x)=\sum_{i=1}^k\alpha^*(x_i)$, $x\in\R^k$.
Combining this bound with (\ref{eq:step1}) leads to
\begin{eqnarray*}
&& \frac{\ell(t)}{\ell'(t)}Z'(t) -Z(t) \log Z(t) \\
&&\qquad\leq
\biggl( \frac{\eta\omega_{\alpha^*} ( {\ell(t)}/{\eta}
)}{1-(K(t) + \eta)C_3} - \frac{\ell(t)^2}{\ell'(t)} \biggr)
\int c^* ( {\nabla P_{1-t} f} )e^{\ell(t) P_{1-t} f}
\,d\mu.
\end{eqnarray*}
Our aim is to choose the various parameters so that to have the
right-hand side of the latter inequality nonpositive.
We will make sure to choose $\ell$ so that $\ell(t)/\eta< 1$; then
by Lemma \ref{lem:tec2} below $K(t)\leq\ell(t)/(1-t)$ and $\omega
_{\alpha^*} ( \frac{\ell(t)}{\eta} )\leq\frac{\ell
^2(t)}{\eta^2}$.
Setting $v=1-C_3\eta$, one has $0< v< 1$,
%
%
\begin{equation}\label{eq:*}
C_3\bigl(K(t)+\eta\bigr)\leq(1-v)\biggl(\frac{\ell(t)}{\eta(1-t)} +1\biggr)
\end{equation}
and
%
%
\begin{eqnarray}\label{inetech}
&&\biggl( \frac{\eta\omega_{\alpha^*} ( {\ell(t)}/{\eta}
)}{1-(K(t) + \eta)C_3} - \frac{\ell(t)^2}{\ell'(t)}
\biggr)\nonumber\\[-8pt]\\[-8pt]
&&\qquad\leq\ell^2(t) \biggl(\frac{1}{\eta v-(1-v){\ell
(t)}/({1-t})}-\frac{1}{\ell'(t)}\biggr).\nonumber
\end{eqnarray}
We choose $\ell(t)=\eta((1-t)^{1-v}-(1-t)), t\in(0,1)$,
so that $\ell(0)=0$ and the right-hand side of (\ref{inetech}) is
equal to zero. Furthermore
$\ell'(t)=\eta(1-\frac{1-v}{(1-t)^v})\geq0, \forall
t\in[0,1-(1-v)^{1/v}]$.
As assumed earlier, $\ell(t)$ is nonnegative and $\ell(t)/\eta< 1$
on $(0,1)$. Let us observe that
\[
\biggl[ \frac{\log Z(t)}{\ell(t)} \biggr]'
=\frac{\ell'(t)}{Z(t)\ell^2(t)}\biggl[\frac{\ell(t)}{\ell
'(t)}Z'(t) -Z(t) \log Z(t) \biggr] .
\]
Let $T=T(v):=1-(1-v)^{1/v}$, since $\ell'(t)>0$ on $(0,T(v))$, the
above inequalities imply that on that interval $[\frac{\log
Z(t)}{\ell(t)}]'\leq0$ provided
$C_3(K(t)+\eta)<1$. By (\ref{eq:*}), this is indeed satisfied for
$t\in[0,T(v)]$.
This gives that the function $t \mapsto\frac{\log Z_t}{\ell(t)}$ is
nonincreasing on $(0,T]$.
Hence, we have
\[
\int e^{\ell(T) P_{T} f} \,d\mu= Z_{T} \leq\exp\biggl( \ell(T)
\lim_{t \to0} \frac{\log Z_t}{\ell(t)} \biggr) = e^{\ell(T) \int
P_1 f \,d\mu} .
\]
In other words, since $P_{T}f \geq f$, then for all bounded continuous
functions $g= \ell(T)f$,
\[
\int e^g \,d\mu\leq e^{\int\tilde P g \,d\mu}
\]
with
\[
\tilde P g (x) = \sup_{y \in\mathbb{R}^k} \{ g(y) - \ell
(T) c(x-y) \} .
\]
According to the Bobkov and G\"otze sup-convolution characterization of
transport inequalities (which for the reader's convenience we quote
below as Theorem \ref{bg}), this implies that $\mu$ verifies
\hyperlink{eqiTcClink}{($\T_c(1/\ell(T))$)}. One has $\ell(T)=\eta v (1-v)^{(1/v)-1}$ and
$C_3\ell(T)=v (1-v)^{1/v}$. Hence, $\mu$ verifies
\hyperlink{eqiTcClink}{($\T_c(K)$)}
with
\[
K = \frac{C_3}{\sup_{v\in(0,1)}v(1-v)^{1/v}}\leq7,7 C_3.
\]
The proof of $(3)\Rightarrow(1)$ is complete.
\end{pf*}
\begin{theorem}[\cite{bobkov-gotze}]\label{bg}
Let $\mu$ be a probability measure on $\mathbb{R}^k$, $\lambda>0$
and $c$ defined as in Theorem \ref{main-result2}.
Then, the following two statements are equivalent:

\begin{enumerate}[(ii)]
\item[(i)] $\mu$ satisfies \hyperlink{eqiTcClink}{($\T_c(1/\lambda)$)};

\item[(ii)] for any bounded function $f\dvtx\mathbb{R}^k \to\mathbb{R}$ it holds
\[
\int e^f \,d\mu\leq
\exp\biggl\{ \int\sup_{y \in\mathbb{R}^k} \{ f(y) - \lambda
c(x-y) \} \biggr\} \,d\mu.
\]
\end{enumerate}
\end{theorem}

Note that Theorem \ref{bg} holds in much more general setting, see
\cite{villani}.

\section{From the restricted logarithmic Sobolev inequality to $\T_{2}$---II:
An alternative proof} \label{sec:alternative}

In this section, we give an alternative proof of the second part $(3)
\Rightarrow(1)$ of Theorem \ref{main-result2}.
The final\vadjust{\goodbreak} result will lead to a worst constant, so we will present our
approach only in the particular case of
the quadratic cost function $c(x)=\frac{1}{2}|x|^2$. More precisely,
we will prove that
(\ref{eqrLSIC}) $\Rightarrow$ \hyperlink{eqiTcClink}{($\T_2(9C)$)}
[leading, for the quadratic cost,
to the implication $(3) \Rightarrow(1)$ of Theorem \ref{main-result2}
with $C_1=9C_3$].
We believe that this alternative approach may lead to other
results in the future and so that it is worth mentioning it.

The strategy is based on the following recent characterization of
Gaussian dimension free concentration by the first author.
\begin{theorem}[\cite{gozlan}]\label{characterization}
A probability measure $\mu$ on $\R^k$ verifies the inequality
\hyperlink{eqiTcClink}{($\T_2(C/2)$)} if and only if there are some $r_o\geq0$ and $b>0$ such that
for all positive integer $n$ and all subset $A$ of $(\R^k
)^n$ with $\mu^n(A)\geq1/2$, the following inequality holds
\[
\mu^n(A+rB_2)\geq1-be^{-(r-r_o)^2/C}\qquad \forall r\geq r_o,
\]
where $B_2$ is the Euclidean unit ball of $(\R^k)^n$.
\end{theorem}

So, in order to get that (\ref{eqrLSIC}) $\Rightarrow$ \hyperlink{eqiTcClink}{($\T_2(9C)$)} it
is enough to
prove that the dimension free Gaussian concentration inequality holds
with $-(r-r_o)^2/(18C)$ in the exponential.

First, let us observe that the restricted logarithmic Sobolev
inequality tensorizes.
\begin{prop}\label{tensorization}
If a probability measure $\mu$ on $\R^k$ verifies (\ref{eqrLSIC})
for some $C>0$, then for all positive integer $n$ the probability $\mu
^n$ verifies (\ref{eqrLSIC}).
\end{prop}
\begin{pf}
If $f\dvtx(\R^k)^n\to\R$ is $K$-semi-convex, then for all
$i\in\{1,\ldots,n\}$ and all $x_1, \ldots, x_{i-1}, x_{i+1},
\ldots, x_n \in\R^k$ the function $f_i\dvtx\R^k\to\R$ defined by
$f_i(x)=f(x_1,\ldots,x_{i-1},x,x_{i+1},\ldots,x_n)$ is $K$-semi-convex.
According to the classical additive property of the entropy functional
(see, e.g., \cite{ane}, Chapter 1)
\[
\ent_{\mu^n}(e^f)\leq\int\sum_{i=1}^n \ent_{\mu}(e^{f_i}) \,d\mu^n.
\]
Applying to each $f_i$ the restricted logarithmic Sobolev inequality
completes the proof.
\end{pf}

The next proposition uses the classical Herbst argument (see,
e.g.,\break
\cite{ledoux}).
\begin{prop}\label{Herbst}
If $\mu$ verifies the restricted logarithmic Sobolev inequality
(\ref{eqrLSIC}) then
for all $f\dvtx\R^k\to\R$ which is $1$-Lipschitz with respect to the
Euclidean norm and $K$-semi-convex with $K\geq0$ one has
\[
\int e^{\lambda(f(x)-\int f \,d\mu)} \,d\mu(x)\leq\exp\biggl(\frac
{2\lambda^2C}{1-\lambda KC}\biggr)\qquad \forall\lambda\in
\bigl(0, 1/(CK)\bigr).
\]
\end{prop}
\begin{pf}
Let us denote $H(\lambda)=\int e^{\lambda f} \,d\mu$, for all $\lambda
\geq0$.
The function $\lambda f$ is $\lambda K$ semi-convex, so if $0\leq
\lambda<1/(CK)$, one can apply the\vadjust{\goodbreak} inequality (\ref{eqrLSIC}) to
the function $\lambda f$.
Doing so yields the inequality
\begin{eqnarray*}
\lambda H'(\lambda)-H(\lambda)\log H(\lambda)
& = &
\ent_\mu( e^{\lambda f} )
\leq
\frac{2C\lambda^2}{(1-\lambda K C)^2}\int|\nabla
f|^2e^{\lambda f} \,d\mu\\
& \leq &
\frac{2C\lambda^2}{(1-\lambda K C)^2}H(\lambda),
\end{eqnarray*}
where the last inequality comes from the fact that $f$ is $1$-Lipschitz.
Consequently, for all $0\leq\lambda<1/(CK)$,
\[
\frac{d}{d\lambda}\biggl(\frac{\log H(\lambda)}{\lambda}
\biggr)\leq\frac{2C}{(1-\lambda K C)^2}.
\]
Observing that $\log H(\lambda)/\lambda\to\int f \,d\mu$ when
$\lambda\to0$ and integrating the differential inequality above gives
the result.
\end{pf}

Now let us show how to approach a given $1$-Lipschitz function by a
$1$-Lipschitz and $K$-semi-convex function.
\begin{prop}\label{approximation}
Let $f\dvtx\R^k\to\R$ be a $1$-Lipschitz function. Define
\[
P_tf(x)=\sup_{y\in\R^k}\biggl\{f(y)-\frac{1}{2t}|x-y|^2\biggr\}\qquad
\forall x\in\R^k, \forall t>0 .
\]
Then:
\begin{enumerate}[(iii)]
\item[(i)] For all $t>0$, $P_t f$ is $1$-Lipschitz.
\item[(ii)] For all $t>0$, $P_t f$ is $1/t$-semi-convex.
\item[(iii)] For all $t>0$ and all $x\in\R^k$, $f(x)\leq P_tf(x)\leq f(x)+
\frac{t}{2}$.
\end{enumerate}
\end{prop}
\begin{pf}
(i)
Write $P_tf(x)=\sup_{z\in\R^k}\{f(x-z)-\frac
{1}{2t}|z|^2\}$. For all $z\in\R^k$, the function $x\mapsto
f(x-z)-\frac{1}{2t}|z|^2$ is $1$-Lipschitz. So $P_tf$ is $1$-Lipschitz
as a supremum of $1$-Lipschitz functions.

(ii)
Expanding $|x-y|^2$ yields $P_t f(x)=\sup_{y\in\R^k}\{
f(y)-\frac{1}{2t}|y|^2+\frac{1}{t}x\cdot y\}-\frac{1}{2t}|x|^2$.
Since a supremum of affine functions is convex, one concludes that
$x\mapsto P_tf(x)+\frac{|x|^2}{2t}$ is convex, which means that $P_tf$
is $1/t$-semi-convex.

(iii)
The inequality $P_tf(x)\geq f(x)$ is immediate. Since $f$ is $1$-Lipschitz,
\begin{eqnarray*}
P_tf(x)-f(x)&=&\sup_{y\in\R^k}\biggl\{f(y)-f(x)-\frac
{1}{2t}|x-y|^2\biggr\}\\ \noalign{\vspace{-2pt}}
&\leq&\sup_{y\in\R^k}\biggl\{|y-x|-\frac{1}{2t}|x-y|^2\biggr\}\\ \noalign{\vspace{-2pt}}
&=&\sup_{r\geq0} \biggl\{r-\frac{r^2}{2t}\biggr\}=\frac{t}{2}
\end{eqnarray*}
\end{pf}

We are now ready to complete the proof.\vadjust{\goodbreak}
\begin{pf*}{Proof of (\ref{eqrLSIC}) $\Rightarrow$ \hyperlink{eqiTcClink}{($\T_2(9C)$)}}
Let $n\geq1$. Consider a $1$-Lipschitz function $g$ on $(\R
^k)^n$ and define
$P_t g(x)=\sup_{y\in(\R^k)^n}\{g(y)-\frac
{1}{2t}|x-y|^2\}$, $t>0$. Thanks to Proposition \ref
{approximation}, the function $P_tg$ is $1$-Lipschitz and
$1/t$-semi-convex, so according to Propositions \ref{tensorization}
and \ref{Herbst}, for all $0\leq\lambda<t/C$, one has
\[
\int e^{\lambda(P_tg(x)-\int P_tg \,d\mu^n)} \,d\mu^n(x)\leq\exp
\biggl(\frac{2\lambda^2C}{1-{\lambda C}/{t}}\biggr).
\]
Moreover, according to point (iii) of Proposition \ref
{approximation}, $P_tg(x)-\int P_tg \,d\mu^n\geq g(x)-\int g \,d\mu
^n-\frac{t}{2}$, for all $x\in(\R^k)^n$. Plugging this
in the inequality above gives
\[
\int e^{\lambda(g(x)-\int g \,d\mu^n)} \,d\mu^n(x)\leq\exp
\biggl(\frac{\lambda t}{2}+\frac{2\lambda^2C}{1-{\lambda
C}/{t}}\biggr).
\]
For a given $\lambda\geq0$, this inequality holds as soon as
$t>C\lambda$. Define $\varphi(t)=\frac{\lambda t}{2}+\frac{2\lambda
^2C}{1-{\lambda C}/{t}}$, $t>0$. It is easy to check that $\varphi
$ attains its minimum value at $t_{\min}=3C\lambda$ (which is
greater than $C\lambda$) and that $\varphi(t_{\min})=9C\lambda^2/2$.
Consequently, we arrive at the following upper
bound on the Laplace transform of $g$:
\[
\int e^{\lambda(g(x)-\int g \,d\mu^n)} \,d\mu^n(x)\leq e^{9C\lambda
^2/2}\qquad \forall\lambda\geq0.
\]
From this, we deduce that every $1$-Lipschitz function $g$ verifies the
following deviation inequality around its mean
\[
\mu^n \biggl(g\geq\int g \,d\mu^n + r\biggr)\leq e^{-r^2/(18C)}\qquad \forall
r\geq0.
\]
Let $r_o$ be any number such that $e^{-r_o^2/(18C)}<1/2$, then denoting
by $m(g)$ any median of $g$, we get $ \int g \,d\mu^n+r_o\geq m(g)$.
Applying this inequality to $-g$, we conclude that $|m(g)-\int g \,d\mu
^n|\leq r_o$. So the following deviation inequality around the median holds
\[
\mu^n\bigl(g\geq m(g)+r\bigr)\leq e^{-(r-r_o)^2/(18C)}\qquad \forall r\geq r_o.
\]
Take $A\subset(\R^k)^n$ with $\mu^n(A)\geq1/2$, and define
$g_A(x)=d_2(x,A)$ where $d_2$ is the usual Euclidean distance. Since
$0$ is a median of $g_A$, the preceding inequality applied to $g_A$ reads
\[
\mu^n(A+rB_2)\geq1-e^{-(r-r_o)^2/(18C)}\qquad \forall r\geq r_o.
\]
According to Theorem \ref{characterization}, this Gaussian dimension
free concentration property implies \hyperlink{eqiTcClink}{($\T_2(9C)$)}.
%
\end{pf*}

\section{Some technical results}\label{sec5}

In this section, we collect some useful results on semi-convex functions.

In the case of differentiable functions, it is easy to rephrase the
definition of semi-convexity, in the following way.
\begin{prop}\label{prop:sc}
Let $c\dvtx\Rk\rightarrow\R^+$ be a differentiable function with
$c(0)=\nabla c(0)=0$. Then,
a differentiable function $f\dvtx\R^k\to\R$ is $K$-semi-convex for the
cost function $c$ if and only if
%
%
\begin{equation}\label{K semi-convex cost c 2}
f(y)\geq f(x)+\nabla f(x)\cdot(y-x)-Kc(y-x)\qquad \forall x,y\in\R^k.
\end{equation}
\end{prop}
\begin{pf}
Suppose that $f$ is $K$-semi-convex; according to the definition, for
all $x,y\in\R^k$ and $\lambda\in[0,1]$, the following holds
\begin{eqnarray*}
f(y)&\geq& f(x)+\frac{f(\lambda x+(1-\lambda)y)-f(x)}{1-\lambda
}\\
&&{}-K\frac{\lambda}{1-\lambda}c\bigl((1-\lambda)(x-y)\bigr)-Kc\bigl(\lambda(y-x)\bigr).
\end{eqnarray*}
Letting $\lambda\to1$ and using $c(0)=\nabla c(0)=0$ one obtains
(\ref{K semi-convex cost c 2}).
Let us prove the converse; according to (\ref{K semi-convex cost c 2}),
\begin{eqnarray*}
f(x)&\geq &f\bigl(\lambda x+(1-\lambda)y\bigr) -(1-\lambda) \nabla f\bigl(\lambda
x+(1-\lambda)y\bigr)\cdot(y-x)\\
&&{}+Kc\bigl((1-\lambda)(y-x)\bigr)
\end{eqnarray*}
and
\[
f(y)\geq f\bigl(\lambda x+(1-\lambda)y\bigr) +\lambda\nabla f\bigl(\lambda
x+(1-\lambda)y\bigr)\cdot(y-x)+Kc\bigl(\lambda(y-x)\bigr).
\]
This gives immediately (\ref{K semi-convex cost c 1}).
\end{pf}
\begin{lem}\label{lemme1}
If $\alpha\dvtx\R\to\R^+$ is a convex symmetric function of class $C^1$
such that $\alpha(0)=\alpha'(0)=0$ and $\alpha'$ is concave on $\R
^+$, then the following inequality holds
%
%
\begin{equation}\label{alpha-ineq}
\alpha(u+v)\leq\alpha(u)+ v\alpha'(u)+4\alpha(v/2)\qquad \forall
u,v\in\R.
\end{equation}
In particular, the function $-c(x) = - \sum_{i=1}^k \alpha(x_i)$,
$x=(x_1,\ldots,x_k)\in\R^k$, is $4$-semi-convex for the cost
$x \mapsto c(x/2)$.
\end{lem}

Note that (\ref{alpha-ineq}) is an equality for $\alpha(t)=t^2$.
\begin{pf*}{Proof of Lemma \ref{lemme1}}
Since $\alpha(v)=\alpha(-v)$, it is enough to prove the inequality
(\ref{alpha-ineq}) for $u\leq0$ and $v\in\R$.
Let us consider the function $G(w):=\alpha(u+w)-\alpha(u)-w\alpha'(u)$.
For $w\geq0$, using the concavity of $\alpha'$ on $\R^+$, either
$u+w\geq0$ and one has
\[
G'(w)=\alpha'(u+w)-\alpha'(u)= \alpha'(u+w)+\alpha'(-u) \leq2
\alpha'(w/2),
\]
or $u+w\leq0$ and one has
\[
G'(w)=\alpha'(-u)-\alpha'(-u-w)\leq\alpha'(w) \leq2 \alpha'(w/2),
\]
since $w\geq0$ and
\begin{eqnarray*}
\frac{\alpha'(w/2)-\alpha'(0)}{w/2}
& \geq &
\frac{\alpha'(w)-\alpha'(0)}{w}\geq\frac{\alpha'(w)-\alpha
'(-u-w)}{2w+u}\\
& \geq &
\frac{\alpha'(-u)-\alpha'(-u-w)}{w}.
\end{eqnarray*}
Similarly, if $w\leq0$, from the convexity of $\alpha'$ on $\R^-$,
$G'(w)\geq\alpha'(w)\geq2\alpha'(w/2)$.
The proof is complete integrating the above inequalities between $0$
and $v$ either for $v\geq0$ or for $v\leq0$.

The second part of the lemma is immediate.
\end{pf*}

The next lemma gives some conditions
on $\alpha$ under which the sup-convolution semi-group $P_t$
transforms functions into semi-convex.
Let us recall that $\omega_\alpha$ is defined by
\[
\omega_\alpha(x) = \sup_{u >0} \frac{\alpha(ux)}{\alpha(u)}
\qquad \forall x\in\R.
\]
\begin{lem} \label{lem:semiconv}
Let $\alpha\dvtx\R\to\R^+$ be a convex symmetric function of class
$C^1$ such that $\alpha(0)=\alpha'(0)=0$ and $\alpha'$ is concave on
$\R^+$.
Let $f\dvtx\R^k\to\R$, $u>0$ and define $g(x)=P_u f(x)=\sup_{y\in\R
^k}\{f(y)-uc((y-x)/u)\}$ with $c(x)=\sum_{i=1}^k \alpha(x_i)$,
$x\in\mathbb{R}^k$. Then
$g$ is $4u \omega_\alpha( \frac{1}{2u} )$-semi-convex
for the cost function $c$.
\end{lem}
\begin{pf}
By Lemma \ref{lemme1}, the function $-c$ is $4$-semi-convex with the
cost function $x \mapsto c(x/2)$.
Consequently, for all $y\in\R^k$, the function $x\mapsto
f(y)-uc((y-x)/u)$ is $4$-semi-convex with the cost function
$x \mapsto uc(x/(2u))$. From the definition (\ref{K semi-convex cost c
1}), we observe that a supremum of $K$-semi-convex functions remains
$K$-semi-convex. Consequently, by definition of $\omega_\alpha$,
we finally get
\begin{eqnarray*}
g(y)
& \geq &
g(x) + \nabla g(x) \cdot(y-x) - 4uc \biggl( \frac{y-x}{2u} \biggr) \\
& \geq &
g(x) + \nabla g(x) \cdot(y-x) - 4u \omega_\alpha\biggl( \frac
{1}{2u} \biggr) c (y-x) .
\end{eqnarray*}
\upqed\end{pf}
\begin{lem}\label{lem:tec2}
Let $\alpha$ be a convex symmetric function of class $C^1$ such that
$\alpha(0)=\alpha'(0)=0$,
$\alpha'$ is concave on $\R^+$. Denote by $\alpha^*$ the conjugate
of $\alpha$.
Then:

\begin{enumerate}[(iii)]
\item[(i)] For any $u \in(0,1)$, $x\in\R$, $\alpha(x/u)\leq\alpha
(x)/u^2$.\vspace*{1pt}

\item[(ii)] For any $u \in(0,1)$, $ \omega_\alpha(1/u ) \leq
1/{u^2}$.

\item[(iii)] For any $u \in(0,1)$, $\omega_{\alpha^*}(u) \leq u^2$.
\end{enumerate}
\end{lem}
\begin{pf}
Point (i). Let $x\geq0$, by concavity of $\alpha'$ on $\R^+$,
$\alpha'(x)\geq u\alpha'(x/u)+(1-u)\alpha'(0)=u\alpha'(x/u)$. The
result follows for $x\geq0$ by integrating between 0 and $x$ and then
for $x\leq0$ by symmetry.
Point (ii) is a direct consequence of point~(i).

Point (iii). Observing that $(\alpha^*)'=(\alpha')^{-1}$, it
follows that $(\alpha^*)'$ is convex on $\R^+$ and $(\alpha
^*)'(0)=\alpha^*(0)=0$. Then the proof is similar to the proof of
point (ii).
\end{pf}

\section{Final remarks}\label{sec6}

In this final section, we state some remarks and extensions on the
topic of this paper.

\subsection{Extension to Riemannian manifolds} \label{Riemannian manifolds}
Otto--Villani theorem holds true on general Riemannian manifolds
\cite{otto-villani}. Furthermore, efforts have been made recently to extend
the Otto--Villani theorem to spaces with poorer structure such as length
spaces \cite{LV07,Balogh09} or general metric spaces \cite{gozlan}.
This section is an attempt to extend our main result to spaces other
than Euclidean spaces. We will focus our attention on the inequality
\hyperlink{eqiTcClink}{($\T_{2}$)} on a Riemannian manifold.

In all what follows, $X$ will be a complete and connected Riemannian
manifold equipped with its geodesic distance $d$:
%
%
\begin{eqnarray}\label{geodesic distance}
d(x,y)=\inf\biggl\{\int_{0}^1 |\dot{\gamma}_{s}| \,ds; \gamma\in
\mathcal{C}^1([0,1],X), \gamma_{0}=x,
\gamma_{1}=y\biggr\}\nonumber\\[-8pt]\\[-8pt]
\eqntext{\forall x,y\in X.}
\end{eqnarray}

A minimizing path $\gamma$ in (\ref{geodesic distance}) is called a
minimal geodesic from $x$ to $y$; in general it is not unique. It is
always possible to consider that minimal geodesics are parametrized in
such a way that
\[
d(\gamma_{s},\gamma_t)=|s-t|d(x,y)\qquad \forall s,t\in[0,1],
\]
and this convention will be in force in all the sequel.



A function $f\dvtx X\to\R$ will be said $K$-semi-convex, $K\geq0$ if for
all $x,y\in X$ and all minimal geodesics $\gamma$ between $x$ and $y$,
the following inequality holds
\[
f(\gamma_s)\leq(1-s)f(x)+sf(y) + s(1-s)\frac{K}{2}d^2(x,y)\qquad
\forall s\in[0,1].
\]
When $f$ is of class $\mathcal{C}^1$ this is equivalent to the
following condition:
%
%
\begin{equation}\label{semiconvex Riem}
f(y)\geq f(x)+\langle\nabla f(x), \dot{\gamma}_0\rangle- \frac
{K}{2}d^2(x,y)\qquad \forall x,y\in X,
\end{equation}
for all minimal geodesics $\gamma$ from $x$ to $y$ (see, e.g.,
\cite{villani}, Proposition 16.2).
If $f$ is semi-convex, then it is locally Lipschitz \cite{villani}.
According to Rademacher's theorem (see, e.g.,
\cite{villani}, Theorem 10.8), $f$ is thus almost everywhere
differentiable. So the
inequality (\ref{semiconvex Riem}) holds for almost all $x\in X$ and
for all $y\in X$.
A function $f$ will be said $K$-semi-concave if $-f$ is $K$-semi-convex.
%
%
\begin{lem}
If $f$ is $K$-semi-convex, then for almost all $x\in X$, the inequality
\[
f(y)\geq f(x)-|\nabla f|(x)d(x,y)-\frac{K}{2}d^2(x,y),
\]
holds for all $y\in X$.
\end{lem}
\begin{pf}
Since the geodesic is constant speed, $|\dot{\gamma}_0|=d(x,y)$.
Applying Cauchy--Schwarz inequality in (\ref{semiconvex Riem}) yields
the desired inequality.
\end{pf}

With this inequality at hand, the proofs of Lemma \ref{lem:easy}
generalizes at once, and we get the following half part of our main result.
\begin{prop}
Suppose that an absolutely continuous probability measure $\mu$ on $X$
verifies the inequality \hyperlink{eqiTcClink}{($\T_2(C)$)}, then it verifies the following
restricted logarithmic Sobolev inequality:
for all $0\leq K<\frac{1}{C}$ and all $K$-semi-convex $f\dvtx X\to\R$,
\[
\ent_{\mu}(e^f)\leq\frac{2C}{(1-KC)^2}
\int|\nabla f|^2 e^f \,d\mu.
\]
\end{prop}

The generalization of the second half part of our main result is more delicate.
We have seen two proofs of the fact that the restricted logarithmic
Sobolev inequality implies \hyperlink{eqiTcClink}{($\T_{2}$)}: one based on the Hamilton--Jacobi
equation and the other based on dimension free concentration.
The common point of these two approaches is that we have used in both
cases the property that the sup-convolution operator $f\mapsto P_t f$
transforms functions into semi-convex functions (see Proposition \ref
{approximation} and Lemma \ref{lem:semiconv}). Let us see how this
property can be extended to Riemannian manifolds.
\begin{lem}\label{P_t semiconvex}
Suppose that there is some constant $S\geq1$, such that the inequality
%
%
\begin{eqnarray}\label{Ohta}
d^2(\gamma_s,y)&\geq&(1-s)d^2(x,y)+sd^2(z,y)\nonumber\\[-8pt]\\[-8pt]
&&{}-s(1-s)S^2d^2(x,z)\qquad
\forall s\in[0,1],\nonumber
\end{eqnarray}
holds for all $x,y,z \in X$, where $\gamma$ is a minimal geodesic
joining $x$ to $z$. This amounts to say that for all $y\in X$, the
function $x\mapsto d^2(x,y)$ is $2S^2$-semi-concave.

Then for all $f\dvtx X\to\R$ and all $u>0$ the function
%
%
\begin{equation}\label{P_t Riem}
x\mapsto P_u f(x)=\sup_{y\in X}\biggl\{ f(y)-\frac{1}{2u}
d^2(x,y)\biggr\}
\end{equation}
is $S^2/ u$-semi-convex.
\end{lem}
\begin{pf}
Under the assumption made on $d^2$, for all $y\in X$, the function
$x\mapsto f(y)-\frac{1}{2u}d^2(x,y)$ is $S^2/u$-semi-convex. Since a
supremum of $S^2/u$ semi-convex functions is $S^2/u$-semi-convex, this
ends the proof.
\end{pf}

Let us make some remarks on condition (\ref{Ohta}). This condition was
first introduced by Ohta in \cite{Ohta09} and Savare in
\cite{Savare07} in their studies of gradient flows in the Wasserstein
space over nonsmooth metric spaces. The condition (\ref{Ohta}) is
related to the Alexandrov curvature of geodesic spaces which
generalizes the notion of sectional curvature in Riemannian geometry.

The first point is a classical consequence of Toponogov's theorem
\cite{Cheeger}.
The second point in the following proposition is due to Ohta
\cite{Ohta09}, Lemma 3.3.
\begin{prop}
Let $X$ be a complete and connected Riemannian manifold.
\begin{enumerate}[(2)]
\item[(1)] The condition (\ref{Ohta}) holds with $S=1$ if and only if the
sectional curvature of $X$ is greater than or equal to $0$ everywhere.
\item[(2)] Suppose that the sectional curvature is greater than or equal to
$\kappa$, where $\kappa\leq0$, then for all $x,y,z \in X$ and every
geodesic $\gamma$ joining $x$ to $z$, one has
%
%
\begin{eqnarray}
d^2(\gamma_s,y)&\geq&(1-s)d^2(x,y)+sd^2(z,y)\nonumber\\[-8pt]\\[-8pt]
&&{}-\Bigl(1+\kappa^2\sup_{t\in[0,1]}d^2(\gamma_t,y)
\Bigr)(1-s)sd^2(x,z).\nonumber
\end{eqnarray}
In particular, if $(X,d)$ is bounded, then (\ref{Ohta}) holds with
\[
S=\bigl(1+\kappa^2\operatorname{diam}(X)^2\bigr)^{1/2}.
\]
\end{enumerate}
\end{prop}

In particular, the case of the Euclidean space, studied in the
preceding sections, corresponds to the case where the sectional
curvature vanishes everywhere.

Now, let us have a look to Hamilton--Jacobi equation. The following
theorem comes from \cite{villani}, Proposition 22.16 and Theorem 22.46.
\begin{theorem}\label{hamilton-generalise}
Let $f$ be a bounded and continuous function on $X$, the function
$(t,x)\mapsto P_t f(x)$ defined by (\ref{P_t Riem}) verifies the
following: for all $t>0$ and $x\in X$,
\[
\lim_{h\to0^+} \frac{P_{t+h}f(x)-P_{t}f(x)}{h}=\frac{|\nabla
^-(-P_{t}f)|^2(x)}{2},
\]
where the metric sub-gradient $|\nabla^- g|$ of a function $g$ is
defined by
\[
|\nabla^- g|(x)=\limsup_{y\to x} \frac{[g(y)-g(x)]_-}{d(y,x)}\qquad
\forall x\in X.
\]
\end{theorem}

Under the condition (\ref{Ohta}), $x\mapsto P_tf(x)$ is semi-convex,
and so differentiable almost everywhere, so for all $t$ and almost all
$x\in X$,
\[
\lim_{h\to0^+} \frac{P_{t+h}f(x)-P_{t}f(x)}{h}=\frac{|\nabla
P_{t}f|^2(x)}{2}.
\]
\begin{theorem}\label{extensionthm}
Suppose that the Riemannian manifold $X$ verifies condition (\ref
{Ohta}) for some $S\geq1$; if an absolutely continuous probability
measure $\mu$ on $X$ verifies the following restricted logarithmic
Sobolev inequality:
for all $0\leq K<\frac{1}{C}$ and all $K$-semi-convex $f\dvtx X\to\R$,
\[
\ent_{\mu}(e^f)\leq\frac{2C}{(1-KC)^2}
\int|\nabla f|^2 e^f \,d\mu,
\]
then it verifies \hyperlink{eqiTcClink}{($\T_2(8CS^2)$)}.
\end{theorem}
\begin{pf}
Setting $C_S=C S^2$, by assumption, for all $KS^2$ semi-convex
functions $f\dvtx X\rightarrow\R$ with $0\leq K<\frac{1}{C_S}$,
\begin{eqnarray*}
\ent_{\mu}(e^f)&\leq&\frac{2C}{(1-KS^2C
)^2} \int|\nabla f|^2 e^f \,d\mu\\
&\leq&\frac{2C_S}{(1-KC_S)^2} \int|\nabla f|^2 e^f
\,d\mu,
\end{eqnarray*}
where the last inequality holds since $S\geq1$. As mentioned in the
\hyperref[intro]{Introduction}, it is still equivalent to \hyperlink{eqrMLSIcClink}{($\mathbf
{rMLSI}(c,C_S)$)} where
$c$ is the quadratic cost function: for all $K\geq0$, $\eta>0$, with
$\eta+K<1/C_S$,
and all $KS^2$ semi-convex functions $f$
%
%
\begin{equation}\label{rmlsi}
\ent_{\mu}(e^f)\leq\frac{\eta}{1-C_S(\eta+ K)} \int
c^*\biggl(\frac{|\nabla f|}{\eta}\biggr) e^f \,d\mu,
\end{equation}
with $c^*(h)= h^2/2$, $h\in\R$.
The end of the proof exactly follows the proof of Theorem \ref
{main-result2} $(3)\Rightarrow(1)$ by replacing $C$ by $C_S$.
There is an additional technical problem due to the right derivatives;
as in the proof of Theorem \ref{main-result2}, we refer to
\cite{LV07,villani} where this difficulty has been circumvented. Therefore,
by Theorem \ref{hamilton-generalise}, we assume that $P_tf$ satisfies
the Hamilton--Jacobi equation $\partial_t P_tf(x) = c^*(|\nabla
P_tf(x)|)$ for all $t>0$ and all $x\in X$. Moreover, by Lemma \ref{P_t
semiconvex} $P_uf$ is $S^2/u$ semi-convex (for the cost
$c(x,y)=d^2(x,y)/2$). Then the continuation of the proof is identical
to the one of Theorem \ref{main-result2} by applying the inequality
(\ref{rmlsi}) to the $K(t)S^2 $ semi-convex function $\ell(t)P_{1-t}f$.
\end{pf}

To conclude this section, let us say that the proof presented in
Section \ref{sec:alternative} can also be adapted to the Riemannian
framework. Essentially,
all we have to do is to adapt the first point of Proposition \ref
{approximation}: the fact that $P_tf$ is $1$-Lipschitz when $f$ is
$1$-Lipschitz. A proof of this can be found in the proof of
\cite{Balogh09}, Theorem 2.5(iv).

\subsection{From transport inequalities to other logarithmic
Sobolev type inequalities} \label{autres log-sob}

Following the ideas of Theorem \ref{th:main1}, we may simply recover
other types of logarithmic Sobolev inequalities. These new forms of
inequalities should be of interest for further developments.
Let $X$ denote a Polish space equipped with the Borel $\sigma
$-algebra. Given Borel functions ${\rmc}\dvtx X\times X\rightarrow
\R$ and $f\dvtx X\rightarrow\R$, define for $\lambda>0$, $x\in X$,
\[
P^\lambda f(x)= \sup_{y\in X} \{ f(y)-\lambda{\rmc}(x,y)
\}.
\]
By definition, one says that a function $ f\dvtx X\rightarrow\R$ is
$K$-semi-concave for the cost ${\rmc}$ if $-f $ is $K$-semi-convex for
the cost ${\rmc}$.
\begin{theorem} \label{th:main1bis}
Let ${\rmc}\dvtx X \times X \rightarrow\R^+$ be a symmetric Borel function.
Let $\mu$ be a probability measure on $X$ satisfying \hyperlink{eqiTcClink}{($\T_{\rmc}(C)$)}
for some
$C>0$. Then for
all $\lambda\in(0,1/C)$, and all function $f\dvtx X\rightarrow\R$,
%
%
\begin{equation}\label{LS1}
\ent_\mu(e^f) \leq\frac{1}{1 - \lambda C} \int
(P^\lambda f - f) \,d\mu\int e^f \,d\mu.
\end{equation}
Assume moreover that ${\rmc}(x,y)=c(x-y)$, $x,y\in\R^k$, where
$c\dvtx\R
^k\rightarrow\R^+ $ is a differentiable symmetric function with
$c(0)=\nabla c(0)=0$. Then
for all $K\geq0 ,\eta> 0$ with $\eta+K<1/C$ and all $K$-semi-concave
differentiable function $f\dvtx\R^k \rightarrow\R$,
%
%
\begin{equation}\label{LS2}
\ent_\mu(e^f) \leq\frac{\eta}{1 - C(\eta+K)} \int
c^*\biggl(\frac{\nabla f}{\eta}\biggr) \,d\mu\int e^f \,d\mu.
\end{equation}
\end{theorem}
\begin{pf}
Following the proof of Theorem \ref{th:main1}, one has for every
probability measure $\pi$ with marginals $\nu_f$ and $\mu$,
\[
H(\nu_{f}| \mu)\leq\iint\bigl(f(x)-f(y)\bigr) \,d\pi(x,y).
\]
From the definition of the sup-convolution function $P^\lambda f$, one has
\[
H(\nu_{f}| \mu)\leq\iint\bigl(P^\lambda f(y) - f(y) \bigr)
\,d\pi(x,y) + \lambda\iint c(y,x) \,d\pi(x,y).
\]
Optimizing over all probability measure $\pi$ and since $\mu$
satisfies (\ref{eqiTcC}), this yields
\[
H(\nu_{f}| \mu)\leq\int\bigl(P^\lambda f(y) - f(y) \bigr)
\,d\mu+ \lambda C H(\nu_{f}| \mu) .
\]
This is exactly the inequality (\ref{LS1}).
Now, if ${\rmc}(x,y)=c(x-y)$, $x,y\in\R^k$, and $f\dvtx\R
^k\rightarrow
\R$ is a $K$-semi-concave differentiable function, then by Lemma \ref
{lem:easy} one has: for all $\eta> 0$,
\[
P^{K+\eta}f- f\leq\eta c^* \biggl(\frac{\nabla f }{\eta} \biggr).
\]
The restricted modified logarithmic Sobolev inequalities (\ref{LS2})
then follows.
\end{pf}

\subsection{\texorpdfstring{On Poincar\'e inequalities}{On Poincare inequalities}} \label{poincare}

Let $c\dvtx\R^k\rightarrow\R$ be a differentiable function such that
$c(0)=\nabla c(0)=0$, with Hessian at point $0$ such that $D^2c(0)>0$
(as symmetric matrices). As for the logarithmic Sobolev inequalities,
it is known that a linearized version of the
transport inequality (\ref{eqiTcC}) is Poincar\'e inequality (see
\cite{maurey,otto-villani,bgl}).

Naturally, (\ref{eqrMLSIcC}) or (\ref{eqICLSIcC}) also
provide Poincar\'e inequality by using basic ideas given in
\cite{maurey} (see also \cite{bgl}). Namely, starting from (\ref{eqICLSIcC}), we apply it with $\varepsilon f$, where $f\dvtx\R
^k\rightarrow\R$ is a smooth function with compact support. The
infimum $\inf_{y \in\mathbb{R}^k} \{ \varepsilon f(y) +
\lambda c ({x-y} ) \} $ is attained at some
$y_\varepsilon$ such that $\varepsilon\nabla f (y_\varepsilon)=
\lambda\nabla c(x-y_\varepsilon)$. Since for $h\in\R^k$, $\nabla
c^* (\nabla c)(h)=h$, one has
\[
x-y_\varepsilon=\nabla c^*\biggl(\frac{\varepsilon\nabla
f(y_\varepsilon)}{\lambda}\biggr)= \frac{\varepsilon}{\lambda}
D^2c^*(0)\cdot\nabla f(x)+ o(\varepsilon).
\]
Therefore, since $D^2c^*(\nabla c(h))\cdot D^2c(h)=I$ and after some
computations, we get the following Taylor expansion
\begin{eqnarray*}
Q^\lambda( \varepsilon f)(x)&=&\varepsilon f(y_\varepsilon)+ \lambda c
( x-y_\varepsilon)\\
&=&\varepsilon f(x)- \frac{\varepsilon^2}{2\lambda} \nabla f(x)^{T}
\cdot D^2c^*(0)\cdot\nabla f(x) + o( \varepsilon^2 ).
\end{eqnarray*}
It is a classical fact that
\[
\ent_\mu( e^{ \ep f} ) = \frac{\varepsilon^2}{2} \Var
_\mu(f) + o( \varepsilon^2 ).
\]
Finally, as $\varepsilon\rightarrow0$, (\ref{eqICLSIcC})
implies: for every $\lambda\in(0,1/C)$,
\[
\Var_\mu(f)\leq\frac{1}{\lambda(1-\lambda C)} \int\nabla f^T\cdot
D^2c^*(0)\cdot\nabla f \,d\mu.
\]
Optimizing over all $\lambda$ yields the following Poincar\'e
inequality for the metric induced by
$D^2c^*(0)$
\[
\Var_\mu(f)\leq4C\int\nabla f^T\cdot D^2c^*(0)\cdot\nabla f \,d\mu.
\]
Denoting by $\|\cdot\|$ the usual operator norm, one also has a
Poincar\'e inequality with respect to the usual Euclidean metric
\[
\Var_\mu(f)\leq4C \|D^2c^*(0)\| \int|\nabla f|^2 \,d\mu.
\]
From the infimum-convolution characterization of transport inequality
(\ref{eqiTcC}) (see Theorem \ref{bg}), a similar proof gives
the same Poincar\'e inequality with the constant $C$ instead of $4C$
(see \cite{maurey}).

Conversely, Bobkov and Ledoux \cite{bobkov-ledoux}, Theorem 3.1,
obtained that Poincar\'e inequality implies a modified logarithmic
Sobolev inequality.
Let $\alpha_{2,1}\dvtx\R\rightarrow\R^+$ and $c_{2,1}\dvtx\R
^k\rightarrow\R^+$ be the cost function defined by
\[
\alpha_{2,1}(h)=\min\bigl( \tfrac{1}{2} h^2, |h| - \tfrac12
\bigr)\qquad \forall h\in\R,
\]
and $c_{2,1}(x)=\sum_{i=1}^k \alpha_{2,1}(x_i), x\in\R^k$. One has
$\alpha_{2,1}^*(h)=h^2/2$ if $|h|\leq1$ and $\alpha
_{2,1}^*(h)=+\infty$ otherwise. Bobkov--Ledoux's result is the following.
\begin{theorem}[\cite{bobkov-ledoux}]\label{bl} Let $\mu$ be a
probability measure on $\R^k$ satisfying the Poincar\'e inequality:
{\renewcommand{\theequation}{$\mathbf{P}(C)$}
\begin{equation}\label{eqPC}
\Var_\mu(f)\leq C \int|\nabla f|^2 \,d\mu,
\end{equation}}

\noindent
for every smooth function $f$ on $\R^k$. Then the following modified
logarithmic Sobo\-lev inequality holds [in short (\ref{eqBLIC})]:
for all $\kappa<2/\sqrt{C}$ and every smooth function $f$,
{\renewcommand{\theequation}{$\mathbf{BLI}(C)$}
\begin{equation}\label{eqBLIC}
\ent_\mu(e^f) \leq C\kappa^2 K(\kappa,C)\int\alpha
_{2,1}^* \biggl(\frac{\nabla f}{\kappa} \biggr) e^f \,d\mu,
\end{equation}}

\noindent
where $K(\kappa,C)= (\frac{2+\kappa\sqrt C}{2-\kappa\sqrt
C})^2 e^{\kappa\sqrt{5C}}$.
\end{theorem}

Applying (\ref{eqBLIC}) to $\varepsilon f$, as $\varepsilon
\rightarrow0$, (\ref{eqBLIC}) yields
$\mathbf{P}( C K(\kappa,C))$ but also
(\ref{eqPC}) since $K(\kappa,C)\rightarrow1$ as $\kappa
\rightarrow0$. Theorem \ref{bl} therefore indicates that $\mathbf
{P}( C)$ and (\ref{eqBLIC}) are exactly equivalent.
Thanks to the Hamilton--Jacobi approach, Bobkov, Gentil and Ledoux
\cite{bgl} obtained that (\ref{eqBLIC}) implies \hyperlink{eqiTcClink}{($\T_{\tilde c^{\kappa
}_{2,1}}(C)$)} for all $\kappa<2/\sqrt C$ where
%
%
\setcounter{equation}{8}
\begin{equation}\label{c-kappa21}
{\tilde c^{\kappa}_{2,1}}(x)=\kappa^2C^2 K(\kappa,C) \alpha
_{2,1}\biggl(\frac{|x|}{\kappa C K(\kappa,C)}\biggr)\qquad \forall
x\in\R^k.
\end{equation}
By linearization and optimization over $\kappa$,
\hyperlink{eqiTcClink}{($\T_{\tilde c^{\kappa}_{2,1}}(C)$)} implies (\ref{eqPC}), and
therefore (\ref{eqBLIC}) is also equivalent to \hyperlink{eqiTcClink}{($\T_{\tilde
c^{\kappa}_{2,1}}(C )$)} for all $\kappa<2/\sqrt C$.

Let\vspace*{1pt} $c^{\kappa}_{2,1}$ denote the cost function defined
similarly as $\tilde c^{\kappa}_{2,1}$ replacing
$\alpha_{2,1}(|\cdot|)$ by $c_{2,1}$ in (\ref{c-kappa21}). One has
$\tilde c^{\kappa}_{2,1}\leq c^{\kappa}_{2,1}$ [this is\vspace*{1pt} a
consequence of the subadditivity of the concave function
$h\rightarrow\alpha_{2,1} (\sqrt h)$]. Therefore,
\hyperlink{eqiTcClink}{($\T_{ c^{\kappa}_{2,1}}(C)$)} implies
\hyperlink{eqiTcClink}{($\T_{\tilde c^{\kappa }_{2,1}}(C)$)}. Consider
now the case of dimension 1, $k=1$, so that $c^{\kappa}_{2,1}=\tilde
c^{\kappa}_{2,1}$.\vspace*{2pt} Theorem \ref{main-result2} indicates
that \hyperlink{eqiTcClink}{($\T_{ c_{2,1}^\kappa}$)} is equivalent, up
to constant, to
\hyperlink{eqrMLSIcClink}{($\mathbf{rMLSI}(c_{2,1}^\kappa)$)}. Actually
\hyperlink{eqrMLSIcClink}{($\mathbf{rMLSI}(c_{2,1}^\kappa)$)} can be
interpreted as $\mathbf {BLI}$ restricted to a class of semi-convex
function for the cost~$c_{2,1}^\kappa$. However, from the discussions
above, \hyperlink{eqrMLSIcClink}{($\mathbf {rMLSI}(c_{2,1}^\kappa)$)}
and $\mathbf{BLI}$ are equivalent up to constant. It would be
interesting to directly recover $\mathbf{BLI}$ from
\hyperlink{eqrMLSIcClink}{($\mathbf{rMLSI}(c_{2,1}^\kappa)$)} or from
\hyperlink{eqiTcClink}{($\T_{ c_{2,1}^\kappa }$)}. The known results
can be summarized by the following diagram for $k=1$:\vspace*{3pt}

%
\noindent\fbox
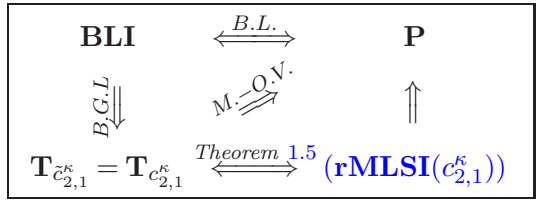
{
\begin{tabular}{@{}cc@{}}
$
\begin{array}{cccc}
\mathbf{BLI} & \stackrel{B.L.}{\Longleftarrow \!\!\Longrightarrow} & \mathbf{P} \\
\begin{sideways} $\!\! \stackrel{B.G.L}{\Longleftarrow}$ \end{sideways}
& \begin{turn}{30} $\stackrel{M.-O.V.}{ \Longrightarrow}$ \end{turn} &
\begin{sideways} $\;\Longrightarrow$ \end{sideways} \\
\T_{\tilde c^{\kappa}_{2,1}} = \T_{ c^{\kappa}_{2,1}} & \stackrel{\fontsize{8.36pt}{10.36pt}
\selectfont{\mbox{\textit{Theorem}~\ref{main-result2}}}}{\Longleftarrow \!\Longrightarrow}
& \mbox{\hyperlink{eqrMLSIcClink}{($\mathbf{rMLSI}(c_{2,1}^\kappa)$)}}
\end{array}
$
\end{tabular}}
\quad
{\small
\begin{tabular}{@{}rl@{}}
where:\hspace*{7pt}\\
\textit{B.L.}: & Bobkov, Ledoux \cite{bobkov-ledoux};\\
\textit{B.G.L.}: & Bobkov, Gentil, Ledoux \cite{bgl};\\
\textit{M.}: & Maurey \cite{maurey};\\
\textit{O.V.}: & Otto, Villani \cite{otto-villani}.
\end{tabular}}

\section*{Acknowledgments}
We warmly thank an anonymous referee for providing constructive
comments and help in improving the contents of this paper.

%

%
\printaddresses

\end{document}